# Design Issues for Generalized Linear Models: A Review

**André I. Khuri, Bhramar Mukherjee, Bikas K. Sinha and Malay Ghosh**

*Abstract.* Generalized linear models (GLMs) have been used quite effectively in the modeling of a mean response under nonstandard conditions, where discrete as well as continuous data distributions can be accommodated. The choice of design for a GLM is a very important task in the development and building of an adequate model. However, one major problem that handicaps the construction of a GLM design is its dependence on the unknown parameters of the fitted model. Several approaches have been proposed in the past 25 years to solve this problem. These approaches, however, have provided only partial solutions that apply in only some special cases, and the problem, in general, remains largely unresolved. The purpose of this article is to focus attention on the aforementioned dependence problem. We provide a survey of various existing techniques dealing with the dependence problem. This survey includes discussions concerning locally optimal designs, sequential designs, Bayesian designs and the quantile dispersion graph approach for comparing designs for GLMs.

## 1. INTRODUCTION

In many experimental situations, the modeling of a response of interest is carried out using regression techniques. The precision of estimating the unknown parameters of a given model depends to a large extent on the design used in the experiment. By design is meant the specification of the levels of the factors (control variables) that influence the response.

The tools needed for the adequate selection of a design and the subsequent fitting and evaluation of the model, using the data generated by the design, have been developed in an area of experimental design known as response surface methodology (RSM). This area was initially developed for the purpose of determining optimum operating conditions in the chemical industry. It is now used in a variety of fields and applications, not only in the physical and engineering sciences, but also in the biological, clinical, and social sciences. The article by Myers, Khuri and Carter (1989) provides a broad review of RSM (see also Myers, 1999). In addition, the three books by Box and Draper (1987), Myers and Montgomery (1995) and Khuri and Cornell (1996) give a comprehensive coverage of the various techniques used in RSM.

*André I. Khuri is Professor, Department of Statistics, University of Florida, Gainesville, Florida 32611, USA e-mail: ufakhuri@stat.ufl.edu. Bhramar Mukherjee is Assistant Professor, Department of Statistics, University of Florida, Gainesville, Florida 32611, USA e-mail: mukherjee@stat.ufl.edu. Bikas K. Sinha is Professor, Division of Theoretical Statistics and Mathematics, Indian Statistical Institute, Kolkata 700 108, India e-mail: sinhabikas@yahoo.com. Malay Ghosh is Distinguished Professor, Department of Statistics, University of Florida, Gainesville, Florida 32611, USA e-mail: ghoshm@stat.ufl.edu.*

Most design methods for RSM models in the present statistical literature were developed around agricultural, industrial and laboratory experiments. These designs are based on the standard general linear model where the responses are assumed to be continuous (quite often, normally distributed) with uncorrelated errors and homogeneous variances. However, clinical or epidemiological data, for example, quite frequently do not satisfy these assumptions. For example, data which consist primarily of human responses tend to be more variable than is expected under the assumption of homogeneous variances. There is much less control over data collected in a clinical setting than over data observed in a laboratory or an industrial setting. Furthermore, most biological data are correlated due to particular genetic relationships. There are also many situations in which clinical experiments tend to yield discrete data. Dose-response experiments are one good example where the responses are binary in most cases. Furthermore, several responses may be observed for the same patient. For example, in addition to the standard binary response of success or failure, some measure of side effects of the treatment may be of importance. Since responses are measured from the same subject, they will be correlated, and hence considering each response to be independent of others may lead to erroneous inferences.

Due to the nature of the data as described above, doing statistical analysis of the data using standard linear models will be inadequate. For such data, generalized linear models (GLMs) would be more appropriate. The latter models have proved to be very effective in several areas of application. For example, in biological assays, reliability and survival analysis and a variety of applied biomedical fields, GLMs have been used for drawing statistical conclusions from acquired data sets. In multicenter clinical trials, estimates of individual hospital treatment effects are obtained by using GLMs (see Lee and Nelder, 2002). In entomology, GLMs are utilized to relate changes in insect behavior to changes in the chemical composition of a plant extract (Hern and Dorn, 2001). Diaz et al. (2002) adopted some GLMs in order to study the spatial pattern of an important tree species. In climatology, GLMs are used to study the basic climatological pattern and trends in daily maximum wind speed in certain regions (see Yan et al., 2002). Also, Jewell and Shiboski (1990) used GLMs to examine the relationship between the risk of HIV (human immunodeficiency virus) infection and the number of contacts with other partners.

In all of the above examples and others, the cornerstone of modeling is the proper choice of designs needed to fit GLMs. The purpose of a design is the determination of the settings of the control variables that result in adequate predictions of the response of interest throughout the experimental region. Hence, GLMs cannot be used effectively unless they are based on efficient designs with desirable properties. Unfortunately, little work has been done in developing such designs. This is mainly due to a serious problem caused by the dependence of a design on the unknown parameters of the fitted GLM.

In this article, we address the aforementioned design dependence problem by providing a survey of various approaches for tackling this problem. The article is organized as follows: Section 2 presents an introduction to GLMs. Section 3 describes criteria for the choice of a GLM design, and introduces the design dependence problem. The various approaches for solving this problem are discussed in Sections 4 (locally optimal designs), 5 (sequential designs), 6 (Bayesian designs) and 7 (quantile dispersion graphs). The article ends with some concluding remarks in Section 8.

## 2. GENERALIZED LINEAR MODELS

As a paradigm for a large class of problems in applied statistics, generalized linear models have proved very effective since their introduction by Nelder and Wedderburn (1972). GLMs are a unified class of regression models for discrete and continuous response variables, and have been used routinely in dealing with observational studies. Many statistical developments in terms of modeling and methodology in the past twenty years may be viewed as special cases of GLMs. Examples include logistic regression for binary responses, linear regression for continuous responses and log-linear models for counts. A classic book on the topic is McCullagh and Nelder (1989). In addition, the more recent books by Lindsey (1997), McCulloch and Searle (2001), Dobson (2002) and Myers, Montgomery and Vining (2002) provide added insights into the application and usefulness of GLMs.

There are three components that define GLMs. These components are:

(i) The elements of a response vector $\mathbf{y}$ are distributed independently according to a certain probability distribution considered to belong to the exponential family, whose probability mass function (or probability density function) is given by

$$(2.1) \quad f(y, \theta, \phi) = \exp\left[\frac{\theta y - b(\theta)}{a(\phi)} + c(y, \phi)\right],$$

where $a(\cdot)$, $b(\cdot)$ and $c(\cdot)$ are known functions; $\theta$ is a canonical parameter and $\phi$ is a dispersion parameter (see McCullagh and Nelder, 1989, page 28).

(ii) A linear regression function, or linear predictor, in $k$ control variables $x_1, x_2, \ldots, x_k$ of the form

$$(2.2) \quad \eta(\mathbf{x}) = \mathbf{f}^T(\mathbf{x})\boldsymbol{\beta},$$

where $\mathbf{f}(\mathbf{x})$ is a known $p$-component vector-valued function of $\mathbf{x} = (x_1, x_2, \ldots, x_k)^T$, $\boldsymbol{\beta}$ is an unknown parameter vector of order $p \times 1$ and $\mathbf{f}^T(\mathbf{x})$ is the transpose of $\mathbf{f}(\mathbf{x})$.

(iii) A link function $g(\mu)$ which relates $\eta$ in (2.2) to the mean response $\mu(\mathbf{x})$ so that $\eta(\mathbf{x}) = g[\mu(\mathbf{x})]$, where $g(\cdot)$ is a monotone differentiable function. When $g$ is the identity function and the response has the normal distribution, we obtain the special class of linear models.

GLMs have several areas of application ranging from medicine to economics, quality control and sample surveys. Applications of the logistic regression model, expanded with the popularity of case-control designs in epidemiology, now provide a basic tool for epidemiologic investigation of chronic diseases. Similar methods have been extensively used in econometrics. Probit and logistic models play a key role in all forms of assay experiments. The log-linear model is the cornerstone of modern approaches to the analysis of contingency table data, and has been found particularly useful for medical and social sciences. Poisson regression models are widely employed to study rates of events such as disease outcomes. The complementary log–log model arises in the study of infectious diseases (e.g., in HIV disease transmission and AIDS as illustrated in Jewell and Shiboski, 1990), and more generally, in the analysis of survival data associated with clinical and longitudinal follow-up studies.

Traditionally, the exponential family model adopted for the study of GLMs deals with a linear function of the response variable involving the unknown parameters of interest. This covers most of the experimental situations arising in practice.

However, some special members, such as the curved exponential family of distributions, are not covered. Thus, GLMs should be further generalized to include such members.

As was pointed out earlier, all known response surface techniques were developed within the framework of linear models under the strong assumptions of normality and equal variances concerning the error distribution. One important area that needs further investigation under the less rigid structure of generalized linear models is the choice of design.

## 3. CHOICE OF DESIGN

By a choice of design we mean the determination of the settings of the control variables that yield an estimated (predicted) response with desirable properties. The mean response, $\mu(\mathbf{x})$, at a point $\mathbf{x}$ in a region of interest, $R$, is given by

$$(3.1) \quad \begin{aligned} \mu(\mathbf{x}) &= h[\mathbf{f}^T(\mathbf{x})\boldsymbol{\beta}] \\ &= h[\eta(\mathbf{x})], \end{aligned}$$

where $\eta(\mathbf{x})$ is the linear predictor in (2.2), and $h$ is the inverse function of the link function $g$. An estimate of $\mu(\mathbf{x})$ is obtained by replacing $\boldsymbol{\beta}$ in (3.1) with $\hat{\boldsymbol{\beta}}$, the maximum likelihood estimate of $\boldsymbol{\beta}$, that is,

$$(3.2) \quad \hat{\mu}(\mathbf{x}) = h[\mathbf{f}^T(\mathbf{x})\hat{\boldsymbol{\beta}}].$$

The variance of $\hat{\mu}(\mathbf{x})$ is approximately given by (see Khuri, 1993, page 198)

$$(3.3) \quad \begin{aligned} &\text{Var}[\hat{\mu}(\mathbf{x})] \\ &= \frac{1}{\phi}\left[\frac{d\mu(\mathbf{x})}{d\eta(\mathbf{x})}\right]^2 \mathbf{f}^T(\mathbf{x})(\mathbf{X}^T\mathbf{W}\mathbf{X})^{-1}\mathbf{f}(\mathbf{x}), \end{aligned}$$

where $\phi$ is the dispersion parameter (determined by the exponential family considered), $\mathbf{X}$ is a matrix whose rows consist of $\mathbf{f}^T(\mathbf{x})$ at the various settings of $\mathbf{x}$ used in a particular design and $\mathbf{W}$ is a diagonal matrix of the form

$$(3.4) \quad \mathbf{W} = \text{diag}(w_1, w_2, \ldots, w_n),$$

where $n$ is the number of experimental runs, and

$$(3.5) \quad w_u = \frac{1}{\phi\sigma_u^2}\left(\frac{d\mu_u}{d\eta_u}\right)^2, \quad u = 1, 2, \ldots, n,$$

where $\sigma_u^2$ is the variance of $y_u$, the response value at the $u$th experimental run, and $\frac{d\mu_u}{d\eta_u}$ denotes the derivative of $\mu(\mathbf{x})$ with respect to $\eta(\mathbf{x})$ evaluated at the setting of $\mathbf{x}$ at the $u$th experimental run ($u = 1, 2, \ldots, n$).

The estimation bias incurred in $\hat{\mu}(\mathbf{x})$ is approximately given by (see Robinson and Khuri, 2003)

$$
\begin{aligned}
(3.6) \quad \text{Bias}[\hat{\mu}(\mathbf{x})] \\
= \mathbf{f}^T(\mathbf{x})(\mathbf{X}^T\mathbf{W}\mathbf{X})^{-1}\mathbf{X}^T\mathbf{W}\boldsymbol{\xi}\frac{d\mu(\mathbf{x})}{d\eta(\mathbf{x})} \\
+ \frac{1}{2\phi}\mathbf{f}^T(\mathbf{x})(\mathbf{X}^T\mathbf{W}\mathbf{X})^{-1}\mathbf{f}(\mathbf{x})\frac{d^2\mu(\mathbf{x})}{d\eta^2(\mathbf{x})},
\end{aligned}
$$

where

$$
\boldsymbol{\xi} = -\frac{1}{2\phi}\mathbf{W}^{-1}\mathbf{Z}_d\mathbf{F}\mathbf{1}_n
$$

and

$$
\mathbf{Z}_d = \text{diag}(z_{11}, z_{22}, \ldots, z_{nn}),
$$

where $z_{uu}$ is the $u$th diagonal element of $\mathbf{Z} = \mathbf{X}(\mathbf{X}^T\mathbf{W}\mathbf{X})^{-1}\mathbf{X}^T$, $\mathbf{F} = \text{diag}(f_{11}, f_{22}, \ldots, f_{nn})$ with

$$
f_{uu} = \frac{1}{\phi\sigma_u^2}\left[\frac{d^2\mu_u}{d\eta_u^2}\right]\left[\frac{d\mu_u}{d\eta_u}\right], \quad u = 1, 2, \ldots, n,
$$

and $\mathbf{1}_n$ is an $n \times 1$ vector of 1's. Here $\frac{d^2\mu_u}{d\eta_u^2}$ denotes the second-order derivative of $\mu(\mathbf{x})$ with respect to $\eta(\mathbf{x})$ evaluated at the $u$th experimental run ($u = 1, 2, \ldots, n$).

A good design is one that minimizes the mean-squared error of $\hat{\mu}(\mathbf{x})$, namely,

$$
\begin{aligned}
(3.7) \quad \text{MSE}[\hat{\mu}(\mathbf{x})] &= E[\hat{\mu}(\mathbf{x}) - \mu(\mathbf{x})]^2 \\
&= \text{Var}[\hat{\mu}(\mathbf{x})] + \{\text{Bias}[\hat{\mu}(\mathbf{x})]\}^2.
\end{aligned}
$$

This is known as the mean-squared error of prediction (MSEP). One major problem in doing this minimization is that the MSEP depends on $\boldsymbol{\beta}$, the parameter vector in the linear predictor in (2.2), which is unknown. This leads us to the so-called design dependence problem. Other design optimality criteria such as $A$-, $D$-, $E$- and $G$-optimality, which are variance-based criteria, suffer also from the same problem.

### 3.1 The Design Dependence Problem

In the foregoing section it has been emphasized that in the context of a GLM, the minimization of the mean-squared error of prediction, or of the variances of the parameter estimators, leading to the so-called optimal designs, depends on the values of unknown parameters. Common approaches to solving this problem include:

(a) The specification of initial values, or best "guesses," of the parameters involved, and the subsequent determination of the so-called locally optimal designs.

(b) The sequential approach which allows the user to obtain updated estimates of the parameters in successive stages, starting with the initial values used in the first stage.

(c) The Bayesian approach, where a prior distribution is assumed on the parameters, which is then incorporated into an appropriate design criterion by integrating it over the prior distribution.

(d) The use of the so-called quantile dispersion graphs approach, which allows the user to compare different designs based on quantile dispersion profiles.

We now provide a review of the basic results that have been developed under the aforementioned approaches.

## 4. LOCALLY OPTIMAL DESIGNS

Binary data under a logistic regression model and Poisson count data are the best known examples to illustrate the implementation of the first approach leading up to a locally optimal design.

### 4.1 Logistic Regression Model

Let us first discuss the study of optimal designs for binary data under a logistic regression model. The key reference is Mathew and Sinha (2001). Other related references include Abdelbasit and Plackett (1983), Minkin (1987), Khan and Yazdi (1988), Wu (1988) and Sitter and Wu (1993). It is postulated that a binary response variable $y$ assumes the values 0 and 1 and the chance mechanism depends on a nonstochastic quantitative covariate $X$ taking values in a specified domain. Specifically, for $X = x, y$ takes the value 1 with probability given by

$$
(4.1) \quad \pi(x) = \frac{1}{1 + \exp(-\alpha - x\beta)},
$$

where $\alpha$ and $\beta$ are unknown parameters with $\beta > 0$. It may be noted that there are other versions of this model studied in the literature. One version refers to the dose-response model which will be discussed in the next sections. The parameters themselves and also some parametric functions such as $\frac{\alpha}{\beta}$ and percentiles of $\pi(x)$ are of interest to the experimenter. The purpose is to suggest continuous (or approximate) optimal designs, that is, optimum dose levels and their relative weights, following the terminology of continuous design theory. See Pukelsheim (1993). It turns out that the solutions to the optimal design problems mentioned above provide optimum values of $\alpha + x\beta$. Hence, in order to implement such designs in practice, good initial estimates, or guess values, of $\alpha$ and $\beta$ are needed. Whereas the earlier researchers established optimality results case by case, Mathew and Sinha (2001) developed a unified approach to tackle the optimality problems by exploiting the property of Loewner-order domination in the comparison of information matrices. However, in some of these studies the technical details depend crucially on the *symmetry* of the transformed factor space. For asymmetric domains, Liski et al. (2002) have initiated some studies in the conventional regression setup. Much yet remains to be done there and also in the context of logistic regression.

With reference to the model (4.1), consider now $s$ distinct dose levels $x_1, x_2, \ldots, x_s$, and suppose we wish to obtain $f_i$ observations on $y$ at dose level $x_i$ $(i = 1, 2, \ldots, s)$. Let $\sum_{i=1}^{s} f_i = n$. For most efficient estimation of $\alpha$ and $\beta$, or some functions thereof, the exact optimal design problem in this context consists of optimally selecting the number $s$ of distinct dose levels, the $x_i$'s (in a given experimental region) and the $f_i$'s, with respect to a given optimality criterion, for a fixed $n$. This is equivalent to the approximate determination of optimum dose levels with respective optimum (relative) weights, denoted by $p_i$ $(i = 1, 2, \ldots, s)$, which sum to 1.

For various application areas, it turns out that the estimation problems that are usually of interest refer to (a) the estimation of $\beta$, or $\alpha/\beta$, or some percentiles of $\pi(x)$ given in (4.1), or (b) the estimation of a pair of parameters such as (i) $\alpha$ and $\beta$, (ii) $\beta$ and $\alpha/\beta$, (iii) $\beta$ and a percentile of $\pi(x)$ and (iv) two percentiles of $\pi(x)$. The approach is to start with the asymptotic variance–covariance matrix of the maximum likelihood estimators of $\alpha$ and $\beta$, and then choose the $x_i$'s and the $p_i$'s optimally by minimizing a suitable function, depending on the nature of the problem at hand and the specific optimality criterion applied. For this reason, we consider the information matrix of the two parameters as a weighted combination of component information matrices based on the $x_i$'s, using the $p_i$'s as weights. Then we argue that the asymptotic variance–covariance matrix of the maximum likelihood estimators of the parameters is just the inverse of the weighted information matrix computed above. The optimality functions to be minimized are different scalar-valued functions of the information matrix. The $D$- and $A$-optimality criteria are well-known examples. Historically, the $D$-optimality criterion has received considerable attention in this context and $A$-optimality has also been considered by some authors; see Sitter and Wu (1993). As was mentioned earlier, the optimum dose levels depend on the unknown parameters $\alpha$ and $\beta$, as is typical in nonlinear settings. In fact, solutions to the optimal design problems mentioned above provide optimum values of $\alpha + \beta x_i$, $i = 1, \ldots, s$. Therefore, while implementing the optimal design in practice, good initial estimates of $\alpha$ and $\beta$ are called for. In spite of this unpleasant feature, it is important to construct the optimal designs in this context; see the arguments in Ford, Torsney and Wu (1992, page 569).

Following the approximate design theory, a design is denoted by $\mathcal{D} = \{(x_i, p_i), i = 1, 2, \ldots, s\}$. Therefore, the information matrix for the joint estimation of $\alpha$ and $\beta$ underlying the design $\mathcal{D}$ is given by

$$
(4.2) \quad \mathbf{I}(\alpha, \beta) = \begin{pmatrix} \sum_{i=1}^{s} p_i \dfrac{\exp(-a_i)}{(1+\exp(-a_i))^2} \\ \sum_{i=1}^{s} p_i x_i \dfrac{\exp(-a_i)}{(1+\exp(-a_i))^2} \\ \sum_{i=1}^{s} p_i x_i \dfrac{\exp(-a_i)}{(1+\exp(-a_i))^2} \\ \sum_{i=1}^{s} p_i x_i^2 \dfrac{\exp(-a_i)}{(1+\exp(-a_i))^2} \end{pmatrix},
$$

where $a_i = \alpha + \beta x_i$, $i = 1, \ldots, s$.

Note that except for the factors $\frac{\exp(-a_i)}{(1+\exp(-a_i))^2}$, the information matrix is identical to the one under the usual linear regression of $y$ on a nonstochastic regressor $x$ under the assumption of homoscedastic error structure. This reminds one of the celebrated de la Garza Phenomenon (de la Garza, 1954) which can be explained as follows. Suppose we consider a $p$th-degree polynomial regression of $y$ on $x$ under homogeneous error structure and we start with an $n$-point design where $n > p + 1$. Then, according to this phenomenon, it is possible to come up with an alternative design with exactly $p + 1$ support points such that the two designs have identical information matrices for the entire set of $p + 1$ parameters. This shows that in the homoscedastic scenario, essentially one can confine attention to the collection of designs supported on exactly $p + 1$ points, if the underlying polynomial regression is of degree $p$. On top of this, it is also possible that a particular $(p + 1)$-point design dominates another $(p + 1)$-point design in the sense that the information matrix (for the entire set of parameters) based on the former design Loewner-dominates that based on the latter. Here, Loewner domination means that the difference of the two information matrices is nonnegative definite. It must be noted that Loewner domination is the best one can hope for, but it is rarely achieved. Pukelsheim (1993) has made some systematic studies of Loewner domination of information matrices. See Liski et al. (2002) for some applications.

In the present setup, however, the form of the information matrix indicates that we are in a linear regression setup, but with a heteroscedastic error structure. Does the de la Garza Phenomenon still hold in this case? More importantly, do we have Loewner domination here? The Mathew and Sinha (2001) article may be regarded as one seeking answers to the above questions. From their study, it turns out that though Loewner domination is not possible, for $D$-optimality the class of two-point designs is *complete* in the entire class of competing designs while for $A$-optimality, it is so in a subclass of competing designs which are *symmetric* in some sense. Here, completeness is in the sense that any competing design outside the class is dominated (with respect to the specific optimality-criterion) by another design within the class.

We shall briefly explain below the salient features of the arguments in Mathew and Sinha (2001). Note first that for the joint estimation of $\alpha$ and $\beta$, or for that matter, for any two nonsingular transforms of them, the $D$-optimality criterion seeks to maximize the determinant of the joint information matrix of $\alpha$ and $\beta$. Routine computations yield an interesting representation for this determinant,

$$
(4.3) \quad
\begin{aligned}
\beta^2 |\mathbf{I}(\alpha, \beta)| &= \left[ \sum_{i=1}^{m} p_i \frac{\exp(-a_i)}{(1 + \exp(-a_i))^2} \right] \\
&\quad \cdot \left[ \sum_{i=1}^{m} p_i a_i^2 \frac{\exp(-a_i)}{(1 + \exp(-a_i))^2} \right] \\
&\quad - \left[ \sum_{i=1}^{m} p_i a_i \frac{\exp(-a_i)}{(1 + \exp(-a_i))^2} \right]^2.
\end{aligned}
$$

Note that for any real number $a$, $\frac{\exp(-a)}{(1+\exp(-a))^2} = \frac{\exp(a)}{(1+\exp(a))^2}$. Moreover, it also turns out that for a given $\mathcal{D} = \{(x_i, p_i), i = 1, 2, \ldots, s\}$, there exists a real number $c > 0$ such that

$$
(4.4) \quad \sum_{i=1}^{s} p_i \frac{\exp(a_i)}{(1 + \exp(a_i))^2} = \frac{\exp(c)}{(1 + \exp(c))^2},
$$

$$
(4.5) \quad \sum_{i=1}^{s} p_i a_i^2 \frac{\exp(a_i)}{(1 + \exp(a_i))^2} \leq c^2 \frac{\exp(c)}{(1 + \exp(c))^2}.
$$

Therefore, for this choice of $c$, the design $\mathcal{D}(c) = [(c, 0.5); (-c, 0.5)]$ provides a larger value of the determinant of the underlying information matrix than that based on $\mathcal{D}$. Hence, the class of two-point *symmetric* designs provides a complete class of $D$-optimal designs. It is now a routine task to determine specific $D$-optimal designs for various parametric functions. Of course, the initial solution is in terms of $c$-optimum (which is independent of $\alpha$ and $\beta$) and then we have to transfer it to optimum dose levels, say $x_0$ and $x_{00}$, by using the relations $c = \alpha + \beta x_0$ and $-c = \alpha + \beta x_{00}$. Thus an initial good guess of $\alpha$ and $\beta$ is called for to evaluate $x_0$ and $x_{00}$.

Again, for $A$-optimality with respect to the parameters $\alpha$ and $\beta$, some algebraic simplification yields the following function (to be minimized) when we restrict to the subclass of *symmetric* designs (i.e., designs involving $a_i$ and $-a_i$ with equal weights for every $i$):

$$
(4.6) \quad
\begin{aligned}
&\mathrm{Var}(\hat{\alpha}) + \mathrm{Var}(\hat{\beta}) \\
&= \left[ \frac{1}{\beta^2} \sum_{i=1}^{m} p_i \frac{\exp(-a_i)}{(1 + \exp(-a_i))^2} [\alpha^2 + \beta^2 + a_i^2] \right] \\
&\quad \cdot \left[ \frac{1}{\beta^2} \sum_{i=1}^{m} p_i \frac{\exp(-a_i)}{(1 + \exp(-a_i))^2} \right. \\
&\quad \left. \cdot \sum_{i=1}^{m} p_i a_i^2 \frac{\exp(-a_i)}{(1 + \exp(-a_i))^2} \right]^{-1} \\
&= \frac{\alpha^2 + \beta^2}{\sum_{i=1}^{m} p_i a_i^2 \exp(-a_i)/(1 + \exp(-a_i))^2} \\
&\quad + \frac{1}{\sum_{i=1}^{m} p_i \exp(-a_i)/(1 + \exp(-a_i))^2}.
\end{aligned}
$$

In view of the existence of the real number $c$ with the properties laid down above, it turns out that the design $\mathcal{D}(c) = [(c, 0.5); (-c, 0.5)]$ once again does better than $\mathcal{D}$ with respect to $A$-optimality, at least in the subclass of symmetric designs so defined. It is a routine task to spell out the nature of a specific $A$-optimal design, which in this case depends on the unknown parameters $\alpha$ and $\beta$ in a twofold manner. First, we have to determine the optimum value for $c$ from given values of the parameters by minimizing the lower bound to (4.6) as a function of $\alpha, \beta$ and $c$. Then we have to evaluate the values of the two recommended dose levels using the relations involving $c$ and $-c$, as above. Note that in the context of $D$-optimal designs, this twofold phenomenon did not arise.

REMARK 4.1.1. It so happens that without the symmetry restriction on the class of competing designs, construction of $A$-optimal designs, in general, is difficult. Numerical computations have revealed that $A$-optimal designs are still *point symmetric* but *not weight symmetric*. An analytical proof of this observation is still lacking. In Mathew and Sinha (2001), $E$-optimal designs have also been studied.

### 4.2 Poisson Count Model

Let us now consider Poisson count data and explain the optimality results following Minkin (1993) and Liski et al. (2002). Here we assume that $y$ follows a Poisson distribution with mean $\mu(x)$, $\mu(x) = c(x) \exp [\theta(x)]$, $\theta(x) = \alpha + \beta x; \beta < 0$. The emphasis in Minkin (1993) was on most efficient estimation of $1/\beta$ by choosing a design in the nonstochastic factor space of $x$ over $[0, \infty)$. Naturally, only a locally optimal design could be characterized, which turned out to be a two-point design with 21.8% mass at 0 and the rest at $2.557\beta$. Thus a good guess for $\beta$ is called for. Liski et al. (2002) developed a unified theory in Minkin's setup for the derivation of a stronger result on completeness of two-point designs, including the point 0. Consequently, it was much easier to re-establish Minkin's result as well as to spell out explicitly the nature of $A$-, $D$-, $E$- and $MV$-optimal designs for simultaneous estimation of the two parameters in the model. It may be noted that Fedorov (1972) gave a complete characterization of $D$-optimal designs in a polynomial regression setup involving several families of heterogeneous variance functions. In an unpublished technical report, Das, Mandal and Sinha (2003) extended the Loewner domination results of Liski et al. (2002) in a linear regression setup involving two specific variance functions. We shall briefly present these results below.

We recall that the setup of Minkin (1993), in a form suitable for our discussion, and as suggested in Liski et al. (2002), is of the form

$$(4.7) \quad \begin{aligned} E(y) &= \alpha + \beta x; \\ V(y) &= v(x)\sigma^2; \\ v(x) &= \exp(x); \quad 0 \leq x < \infty. \end{aligned}$$

As usual, we assume $\sigma^2 = 1$ and confine attention to approximate design theory. Thus, as before, to start with we have an $s$-point design $\mathcal{D}_s = [(x_i, p_i); 0 \leq x_1 < x_2 < \cdots < x_s; \sum_i p_i = 1]$ for $s \geq 2$. It is a routine task to write down the form of the information matrix for the parameters $\alpha, \beta$ underlying the design $\mathcal{D}_s$. We call it $\mathbf{I}_{\mathcal{D}_s}$. Liski et al. (2002) established that given $\mathcal{D}_s$, one can construct a two-point design $\mathcal{D}_2^*$ whose information matrix, say $\mathbf{I}_{\mathcal{D}_2^*}$, dominates $\mathbf{I}_{\mathcal{D}_s}$. Das, Mandal and Sinha (2003) generalized this result when the above form of $v(x)$ is changed to

$$(4.8) \quad \begin{aligned} &\text{(i)} \quad v(x) = k^x, \quad k \geq 1; \\ &\text{(ii)} \quad v(x) = (1+x)^{(1+\gamma)/2}, \quad \gamma \geq -1. \end{aligned}$$

Following Das, Mandal and Sinha (2003), we start with a general variance function $v(x)$ subject to $v(0) = 1$ and $v(x)$ increasing in $x$ over 0 to $\infty$. Next, we start with a two-point design of the form $[(a, p); (b, q)]$ where $0 < a < b < \infty$ and $0 < p, q = 1 - p < 1$. Then, we ask for the sort of variance functions for which this two-point design can be Loewner-dominated by another two-point design, including the point 0. Define in this context two other related functions, $\psi(x) = 1/v(x)$ and $\phi(x) = (v(x) - 1)/x$ for every $x > 0$, $\phi(0)$ being the limit of $\phi(x)$ as $x$ tends to 0. It follows that whenever these two functions satisfy the following conditions, it is possible to achieve this fit:

(i) $\phi(0) = 1$;
(ii) $\phi(x)$ is increasing in $x$;
(iii) for some $s$ and $c$, $0 < s < 1$ and $c > 0$, for which

$$1 - s = [p(1 - \psi(a)) + q(1 - \psi(b))]/[1 - \psi(c)],$$

$\psi(c)$ satisfies the inequality $\psi(c) < p\psi(a) + q\psi(b)$.

Das, Mandal and Sinha (2003) demonstrate that for both forms of $v(x)$ as in (4.8), the above conditions are satisfied. They then continue to argue that this result on Loewner domination (by a two-point design, including the point 0) holds even when one starts with an $s$-point design for $s > 2$. This greatly simplifies the search for specific optimal designs under different variance structures covered by the above two forms. Minkin's result follows as a special case of (i) in the above model. The details are reported in Das, Mandal and Sinha (2003).

## 5. SEQUENTIAL DESIGN

In the previous section, initial values of the parameters are used as best "guesses" to determine a locally optimal design. Response values can then be obtained on the basis of the generated design. In the sequential approach, experimentation does not stop at this initial stage. Instead, using the information thus obtained, updated estimates of the parameters are developed and then used to determine additional design points in subsequent stages. This process continues until convergence is achieved with respect to a certain optimality criterion, for example, $D$-optimality. The implementation of such a strategy is feasible provided that the response values in a given stage can be obtained in a short time, as in sensitivity testing. Sequential designs for GLMs were proposed by Wu (1985), Sitter and Forbes (1997) and Sitter and Wu (1999), among others.

The theory of optimal designs, as noted earlier, involves selection of design points with the goal of minimizing some objective function which can often be interpreted as an expected loss. In a sequential framework, with the arrival of each data point, one needs to make a decision of whether or not to pay for the cost of additional data or else to stop sampling and make a decision. Also, if additional observations become necessary, when there is an option of selecting samples from more than one population, the choice of a suitable sampling rule is equally important. The latter, on many occasions, amounts to the selection of suitable design points. Thus the issue of optimal designs goes hand in hand with sequential analysis, and together they constitute an area of what has become known as sequential design of experiments.

Research on sequential designs can be classified into multiple categories depending on the objective of the researcher. We review here primarily a broad area commonly referred to as "Stochastic Approximation," which was initiated by Robbins and Monro (1951), and was subsequently extended by numerous authors. Here, the problem is one of sequential selection of design points according to some optimality criterion. We also discuss briefly some work on multistage designs.

In Section 5.1, we begin with the stochastic approximation procedure as described in Robbins and Monro (1951). We then consider several extensions and modifications of this pioneering work, and discuss asymptotic properties of the proposed methods. The Robbins–Monro article and much of the subsequent work do not prescribe any stopping rule in this sequential experimentation. There are a few exceptions, and we point out one such result.

Section 5.2 discusses application of the Robbins–Monro method and its extensions for estimating the percentiles of the quantal response curve, in particular, estimation of the median effective dose, pop-ularly known as ED50. We highlight in this context the work of Wu (1985). We also discuss briefly the Bayesian stopping rule as proposed by Freeman (1970). Finally, in Section 5.3, we provide a very short account of multistage designs.

In general, the core of sequential analysis involves sample size determination. This itself is an "optimal" design problem as stopping rules, in general, are motivated by some optimality criteria. A very succinct account of sequential design of experiments, motivated by several important statistical criteria, appeared in Chernoff's (1972) classic monograph.

## 5.1 Stochastic Approximation

We begin with a description of the stochastic approximation procedure of Robbins and Monro (RM) (1951). Let $y$ be a random variable such that conditional on $x$, $y$ has a distribution function (df) $H(y|x)$ with mean $E(y|x) = M(x)$ and variance $V(y|x) = \sigma^2(x)$. The function $M(x)$ is unknown to the experimenter, but is assumed to be strictly increasing so that the equation $M(x) = \alpha$ has a unique root, say $\zeta$. It is desired to estimate $\zeta$ by making successive observations on $y$, say, $y_1, y_2, \ldots$ at levels $x_1, x_2, \ldots$. The basic problem is the selection of the design points $x_1, x_2, \ldots$.

RM address this problem as follows: consider a nonstationary Markov chain with an arbitrary initial value $x_1$, and then define recursively

$$(5.1) \quad x_{n+1} = x_n - a_n(y_n - \alpha), \quad n = 1, 2, \ldots,$$

where $y_n$ has the conditional distribution function $H(y|x_n)$, and $\{a_n\}$ is a sequence of positive constants satisfying

$$0 < \sum_{n=1}^{\infty} a_n^2 < \infty.$$

A finite sample justification of the RM procedure is given by Wu (1985) when $y$ and $x$ are related by a simple linear regression model. In particular, let $y_i = \beta_0 + \beta_1 x_i + e_i$, where $e_i$ are independent and identically distributed with mean 0. Then the parameter of interest is $\zeta = -\beta_0/\beta_1$, the solution of $\beta_0 + \beta_1 x = 0$. If $\beta_1$ is known, $\hat{\beta}_{0n} = \bar{y}_n - \beta_1 \bar{x}_n$ is the least-squares estimator of $\beta_0$ based on $\{(x_i, y_i); i = 1, \ldots, n\}$. Thus a natural choice of $x_{n+1}$ is given by $x_{n+1} = -\hat{\beta}_{0n}/\beta_1 = \bar{x}_n - \beta_1^{-1} \bar{y}_n$. It is shown in Lai and Robbins (1979) that $x_{n+1} = \bar{x}_n - \beta_1^{-1} \bar{y}_n$ for all $n$ is equivalent to $x_{n+1} = x_n - (\beta_1^{-1}/n)y_n$ for all $n$. This is an RM procedure with $a_n = \beta_1^{-1}/n$ and $\alpha = 0$.

RM investigated asymptotic properties of their procedure. Let $b_n = E(x_n - \zeta)^2$. RM showed that for bounded $y$ and strictly increasing $M(x)$, if $a_n = O(n^{-1})$, then $b_n \to 0$ as $n \to \infty$.

A more general result is proved later in Robbins and Siegmund (1971). Rather than the boundedness of $y$ and the monotonicity of $M(x)$, one assumes boundedness of $\sigma(x) + |M(x)|$ by a linear function of $|x|$. These conditions neither imply nor are implied by the original RM conditions. In addition, Robbins and Siegmund (1971) have another assumption which essentially implies that $M(x)$ and $x - \zeta$ are of the same sign. Finally, they assume that

$$(5.2) \qquad 0 < \sum_{n=1}^{\infty} a_n = \infty, \quad 0 < \sum_{n=1}^{\infty} a_n^2 < \infty,$$

which is strictly weaker than the condition $a_n = O(n^{-1})$ as required in RM.

Fabian (1968) proved the asymptotic normality of $x_n$ under additional assumptions. His conditions are similar to those of Robbins and Siegmund, but (5.2) is replaced by the stronger condition

$$(5.3) \qquad n^\gamma a_n \to a(> 0)$$

as $n \to \infty$ for some $\gamma \in (1/2, 1]$. Also, he needed a Lindeberg-type condition for the second moment of the conditional distribution of $y$ given $x$ and also required $M$ to have a positive derivative $M'(\zeta)$ at $\zeta$ and $M'(\zeta) > (2a)^{-1}$ when $\gamma = 1$. The variance of the asymptotic distribution depended on $\sigma(\zeta)$, $a$ and $M'(\zeta)$. The asymptotic variance was given by $a\sigma^2(\zeta)/(2M'(\zeta))$ if $\gamma \in (1/2, 1)$ and $a^2\sigma^2(\zeta)/[2aM'(\zeta) - 1]$ if $\gamma = 1$.

REMARK 5.1.1. The best rate of convergence is achieved when $\gamma = 1$, and then the minimum asymptotic variance is given by $\sigma^2(\zeta)/[M'(\zeta)]^2$. This is attained when $a = [M'(\zeta)]^{-1}$.

In view of Remark 5.1.1, an optimal asymptotic choice of $\{a_n\}$ is given by $a_n = (n + n_0)^{-1} d_n$, where $d_n$ is any consistent estimator of $M'(\zeta)$ and $n_0$ is positive. One can interpret $n_0$ as a "tuning parameter" which does not affect the asymptotics at all, but can play a major role when the sample size is small. Also, Fabian (1983) has shown that when the conditional probability density function of $y$ given $x$ is $N(M(x), \sigma^2(x))$, then the RM procedure is locally asymptotically minimax. Moreover, Abdelhamid (1973) and Anbar (1973) pointed out independently that even for nonnormal conditional probability density functions, the RM process can be made asymptotically optimal by a suitable transformation of observations.

The next thing is to outline procedures which guarantee consistent estimation of $M'(\zeta)$. Venter (1967) addressed this by taking observations in pairs at $x_i \pm c_i$ for a suitable positive sequence of constants $\{c_i\}, i \geq 1$. We denote by $y_{i1}$ and $y_{i2}$ the corresponding responses so that

$$(5.4) \qquad \begin{aligned} E(y_{i1}|x_i) &= M(x_i - c_i); \\ E(y_{i2}|x_i) &= M(x_i + c_i). \end{aligned}$$

Let $z_i = (y_{i2} - y_{i1})/(2c_i)$ and $y_i = (y_{i2} + y_{i1})/2$. If $x_i \xrightarrow{P} \zeta$, then estimate $M'(\zeta)$ by $d_n = \bar{z}_n = \frac{1}{n}\sum_{i=1}^{n} z_i$. Define recursively the design points

$$(5.5) \qquad x_{n+1} = x_n - (n\hat{d}_n)^{-1}(y_n - \alpha), \quad n = 1, 2, \ldots,$$

where $\hat{d}_n$ is a truncated version of $d_n$ ensuring that $x_n \xrightarrow{P} \zeta$. Fabian (1968) and later Lai and Robbins (1979) contain relevant asymptotic results.

In spite of all the asymptotic niceties, the RM procedure can be seriously deficient for finite samples. This was demonstrated in the simulation studies of Wu (1985) and Frees and Ruppert (1990). First, the initial choice of $x_1$ heavily influences the recursive algorithm. If $x_1$ is far away from $\zeta$, then the convergence of $x_n$ to $\zeta$ may demand prohibitively large $n$. Again, if all the $x_i$'s are closely clustered around $\zeta$, no estimator of $M'(\zeta)$ can be very accurate, and imprecise estimation of $M'(\zeta)$ leads to imprecise estimation of $\zeta$ due to the proposed recursive algorithm, at least for small samples.

The original RM paper, and much of the subsequent literature, address sequential design problems without specifying any stopping rule. However, a simple stopping rule can be proposed based on asymptotic considerations. In particular, in the setup of Venter (1967), $\hat{d}_n$ consistently estimates $M'(\zeta)$. Defining $\hat{\sigma}_n^2 = (2n)^{-1}\sum_{i=1}^{n}(y_{i1}^2 + y_{i2}^2)$, it can be shown that $\hat{\sigma}_n^2$ is a consistent estimator of $\sigma^2(\zeta)$. Now invoking the asymptotic normality of $x_n$, a large sample $100(1-\gamma)\%$ confidence interval for $\zeta$ is given by $x_n \pm \Phi^{-1}(1 - \gamma/2)\hat{\sigma}_n/\{\hat{d}_n(2n)^{1/2}\}$, where $\Phi$ is the cumulative distribution function of the $N(0,1)$ variable. If one decides to construct a confidence interval of length not exceeding $l$, then one can stop at the first $n$ for which the width of the above interval is less than or equal to $l$. Sielken (1973) proposed this stopping rule and proved that as $l \to 0$ (so that $n \to \infty$), the coverage probability converges to $1 - \gamma$. Alternative stopping rules have been proposed by Stroup and Braun (1982) and Wei (1985) when the responses are normal. Ruppert et al. (1984) proposed a stopping rule for Monte Carlo optimization which is so designed that $M(x_n)$ is within a specified percentage of $M(\zeta)$.

Kiefer and Wolfowitz (1952) considered a related problem of locating the point where the regression function is maximized or minimized instead of finding the root of a regression equation as in RM. This amounts to solving $M'(\zeta) = 0$. The resulting procedure is similar to that of RM, but $y_n - \alpha$ in (5.1) is replaced by an estimate of $M'(\zeta)$. We omit the details.

One of the major applications of the RM stochastic approximation procedure is the estimation of the percentiles of the quantal response curve. We shall specifically discuss this problem in the next section.

### 5.2 Quantal Response Curves

Consider now the special case when $y$ is a binary variable with

$$(5.6) \quad E(y|x) = \pi(x) = [1 + \exp\{-\theta(x - \mu)\}]^{-1}.$$

The independent variable $x$ is the dose level, while the binary response $y$ is a quantal variable. Such a model is called a "dose-response" model. The parameter $\mu$ is the 50% response dose, or ED50. It is easy to recognize this as a logistic model with location parameter $\mu$ and scale parameter $\theta$ and as a simple reparameterization of the model given in (4.1). Also, here $\pi(x)$ is the same as $M(x)$. Equating $\pi(x)$ to $\alpha$ leads to the solution $\zeta = \mu + \theta^{-1} \log(\frac{\alpha}{1-\alpha})$. Hence, efficient estimation of $\zeta$ depends on efficient estimation of $\mu$ and $\theta$.

The direct RM procedure will continue to generate the $x$-sequence with an arbitrary initial value $x_1$, and then define $x_{n+1} = x_n - a_n(y_n - \alpha)$, where one may take $a_n$ as proportional to $n^{-1}$. However, as recognized by many, including RM, for specific distributions there may be a more efficient way of generating the $x$ values. Indeed, simulations by Wu (1985), and subsequently by Frees and Ruppert (1990), demonstrate the poor performance of the RM procedure for small $n$ in this example.

Wu's procedure begins with generating some initial values $x_1, \ldots, x_m$. Then one obtains the MLE's $\hat{\mu}_m$ and $\hat{\theta}_m$ of $\mu$ and $\theta$, respectively, based on $\{(x_i, y_i); i = 1, \ldots, m\}$ by solving the likelihood equations (i) $\sum_{i=1}^m y_i = \sum_{i=1}^m \pi(x_i)$ and (ii) $\sum_{i=1}^m x_i y_i =$

$\sum_{i=1}^m x_i \pi(x_i)$. Let $x_{m+1} = \hat{\mu}_m + \hat{\theta}_m^{-1} \log(\frac{\alpha}{1-\alpha})$. Update the MLE's of $\mu$ and $\theta$ by $\hat{\mu}_{m+1}$ and $\hat{\theta}_{m+1}$, respectively, based on $\{(x_i, y_i); i = 1, \ldots, m+1\}$. Now let $x_{m+2} = \hat{\mu}_{m+1} + \hat{\theta}_{m+1}^{-1} \log(\frac{\alpha}{1-\alpha})$, and continue in this manner. The $\{x_n\}$ generated in this way meet the consistency and asymptotic optimality properties mentioned in the previous section. Wu (1985) also suggested modification of the likelihood equations (i) and (ii) by multiplying both sides of the $i$th component by a factor $w(|x_i - x_m|)$, where $w$ is a certain weight function. This can partially overcome vulnerability of the logits at extreme tails.

As mentioned earlier, the case $\alpha = 1/2$ is of special interest since then $\zeta = \mu = \text{ED50}$. In this case, an alternative procedure known as the "up-and-down" procedure for generating the design points was proposed by Dixon and Mood (1948). Specifically, let

$$x_{n+1} = \begin{cases} x_n + \Delta, & \text{if } y_n = 0, \\ x_n - \Delta, & \text{if } y_n = 1, \end{cases}$$

where $\Delta$ is somewhat arbitrary. Once again, performance of this procedure depends very much on the choice of a good guess for $x_1$ and $\Delta$. Unless $\Delta$ is made adaptive, the large sample performance of $x_n$ cannot be studied. Wetherill (1963) discusses some modifications of this method.

One important issue that has not been discussed so far is the choice of the stopping rule. This requires explicit consideration of payoff between the cost of further observation and that of less accurate estimation. The study was initiated by Marks (1962) who considered the problem as one of Bayesian sequential design with known $\theta$ and a two-point prior for $\mu$. Later, Freeman (1970) considered the problem, once again for known $\theta$, but with a conjugate prior for $\mu$. In addition, he considered special cases of one, two or three dose levels. For one dose level, he introduced the prior

$$p(\mu|r_0, n_0) \propto \frac{\exp(r_0 \theta(\zeta - \mu))}{[1 + \exp(\theta(\zeta - \mu))]^{n_0}}.$$

This amounts to a Beta$(r_0, n_0 - r_0)$ prior for the response probability $\alpha$ at dose level $\zeta$, with prior parameters $0 < r_0 < n_0$.

With squared error loss plus cost and a uniform prior, that is, $r_0 = 1$ and $n_0 = 2$, Freeman (1970) considered the usual backward induction argument to set up the necessary equations for finding the Bayes stopping rule. The stopping rule cannot be found analytically, but Freeman provided the numerical algorithm for finding the solution. He also considered situations with two or three possible doses, but did not provide a general algorithm. Thus, what is needed here is an approximation of the "optimal" stopping rule. A Bayesian approach using dynamic programming is hard to implement in its full generality. It appears that a suitable approximation of the Bayes stopping rule which retains at least the asymptotic optimality of the regular Bayes stopping rule is the right approach toward solving this problem.

### 5.3 Multistage Designs

In most practical situations, it is more realistic to adopt a multistage design rather than a sequential design, since continuous updating of the inferential procedure with each new observation may not be very feasible. We discuss one such application as considered in Storer (1989) in the context of Phase I clinical trials. These trials are intended to estimate the maximum tolerable dose (MTD) of a new drug. Although a strict quantitative definition of MTD does not exist in clinical trials, very often the 33rd percentile of the tolerance distribution is taken to define MTD.

From a clinician's perspective, an optimal design is one in which the MTD is defined by the dose at which the trial stops. Storer (1989, 1990) argues that the design problem should be viewed instead as an efficient way of generating samples, wherein the design and analysis are robust to the vagaries of patient treatment in a clinical setting. He first introduces four single-stage designs, and then proposes two two-stage designs by combining some of these single-stage designs. The details are available in the two cited papers. Storer implemented his proposal for MTD estimation in a dose-response setting with three logistic curves.

## 6. BAYESIAN DESIGNS

Application of Bayesian design theory to generalized linear models is a promising route to avoid the design dependence problem. One approach, as discussed in Section 3, is to design an experiment for a fixed best guess of the parameters leading to a "locally optimal" (Chernoff, 1953) design. Locally optimal designs for nonlinear models were first suggested in the seminal paper by Box and Lucas (1959). As suggested in Box and Lucas (1959), another natural approach to solve this problem is to express the uncertainty in the parameters through a prior distribution on the parameters. The Bayesian design problem for normal linear models has been discussed in Owen (1970), Brooks (1972, 1974, 1976, 1977), Chaloner (1984), Pilz (1991) and DasGupta (1996). Chaloner and Verdinelli (1995) present an excellent overview of Bayesian design ideas and their applications. Atkinson and Haines (1996) discuss local and Bayesian designs specifically for nonlinear and generalized linear models.

### 6.1 Bayesian Design Criteria

The Bayesian design criteria are often integrated versions of classical design optimality criteria where the integration is carried out with respect to the prior distribution on $\boldsymbol{\beta}$ [$\boldsymbol{\beta}$ is as in (2.2)]. Most of the Bayesian criterion functions are based on normal approximations to the posterior distribution of the vector of parameters $\boldsymbol{\beta}$, as computations involving the exact posterior distribution are often intractable. Several such approximations to the posterior are available (Berger, 1985) and involve either the observed or the expected Fisher information matrix. The most common form of such normal approximation states that under standard regularity conditions, the posterior distribution of $\boldsymbol{\beta}$ is $N(\widehat{\boldsymbol{\beta}}, [n\mathbf{I}(\widehat{\boldsymbol{\beta}}, \xi)]^{-1})$, where $\widehat{\boldsymbol{\beta}}$ is the maximum likelihood estimate (MLE) of $\boldsymbol{\beta}$, and for any given design measure $\xi$, the expected Fisher information matrix is denoted by $\mathbf{I}(\boldsymbol{\beta}, \xi)$. We shall assume, as usual, that the design measure $\xi$ puts relative weights $(p_1, p_2, \ldots, p_k)$ at $k$ distinct points $(x_1, \ldots, x_k)$, respectively, with $\sum_{i=1}^{k} p_i = 1$. The prior distribution for $\boldsymbol{\beta}$ is used as the predictive distribution of $\widehat{\boldsymbol{\beta}}$ and consequently the Bayesian optimality criteria can be viewed as approximations to the expected posterior utility functions under the prior distribution $p(\boldsymbol{\beta})$. The first criterion which is an analogue to the $D$-optimality criterion in classical optimality theory is given by

$$(6.1) \qquad \phi_1(\xi) = E_{\boldsymbol{\beta}}[\log \det \mathbf{I}(\boldsymbol{\beta}, \xi)].$$

Maximizing this function is equivalent to approximately maximizing the expected increase in the Shannon information or maximizing the expected Kullback–Leibler distance between the posterior and prior distributions (Lindley, 1956; DeGroot, 1986; Bernardo, 1979).

The next criterion is of interest when the only quantity to be estimated is a function of $\boldsymbol{\beta}$, say $h(\boldsymbol{\beta})$. In such situations, the approximate asymptotic variance of the MLE of $h(\boldsymbol{\beta})$ is

$$\mathbf{c}(\boldsymbol{\beta})^T [\mathbf{I}(\boldsymbol{\beta}, \xi)]^{-1} \mathbf{c}(\boldsymbol{\beta}),$$

where the $i$th element of $\mathbf{c}(\boldsymbol{\beta})$ is $c_i(\boldsymbol{\beta}) = \partial h(\boldsymbol{\beta})/\partial \beta_i$. The Bayesian $c$-optimality criterion approximates the posterior expected utility (assuming a squared error loss) under the prior $p(\boldsymbol{\beta})$ as

$$(6.2) \quad \phi_2(\xi) = -E_{\boldsymbol{\beta}}\{\mathbf{c}(\boldsymbol{\beta})^T[\mathbf{I}(\boldsymbol{\beta}, \xi)]^{-1}\mathbf{c}(\boldsymbol{\beta})\}.$$

As in the generalization from $c$- to $A$-optimality in classical optimality theory, if one is interested in estimating several functions of $\boldsymbol{\beta}$ with possibly different weights attached to them and if $\mathbf{B}(\boldsymbol{\beta})$ is the weighted average of the individual matrices of the form $\mathbf{c}(\boldsymbol{\beta})\mathbf{c}(\boldsymbol{\beta})^T$, then the criterion to be maximized is

$$(6.3) \quad \phi_3(\xi) = -E_{\boldsymbol{\beta}}\{\operatorname{tr} \mathbf{B}(\boldsymbol{\beta})[\mathbf{I}(\boldsymbol{\beta}, \xi)]^{-1}\}.$$

The functions $\mathbf{c}(\boldsymbol{\beta})$, $\mathbf{B}(\boldsymbol{\beta})$ may not depend on $\boldsymbol{\beta}$ if one is considering only linear functions of $\boldsymbol{\beta}$.

Since the interpretation of the Bayesian alphabetic optimality criteria, as approximations to expected utility, is based on normal approximations to the posterior, Clyde and Chaloner (1996, 2002) suggest several approaches to verify normality through imposing constraints and discuss how to attain such multiple design objectives in this context. Other design criteria which can be related to a Bayesian perspective appear in Tsutakawa (1972, 1980), Zacks (1977) and Pronzato and Walter (1987). How well the Bayesian criteria actually approximate the expected utility in small samples is not very well known. Some illustrations are presented in Atkinson et al. (1993), Clyde (1993a) and Sun, Tsutakawa and Lu (1996). Dawid and Sebastiani (1999) attempt to connect this type of criterion-based Bayesian design to a purely decision-theoretic utility-based approach to Bayesian experimental design.

Müller and Parmigiani (1995) and Müller (1999) suggest estimating the exact posterior utility through Markov chain Monte Carlo (MCMC) methods instead of using analytical approximations to the posterior distribution. They embed the integration and maximization of the posterior utility function by curve fitting of Monte Carlo samples. This is done by simulating draws from the joint parameter/sample space and evaluating the observed utilities and then fitting a smooth surface through the simulated points. The fitted surface acts as an estimate to the expected utility surface and the optimal design can then be found deterministically by studying the extrema of this surface. Parmigiani and Müller (1995), Clyde, Müller and Parmigiani (1995) and Palmer and Müller (1998) contain applications of these simulation-based stochastic optimization techniques.

## 6.2 Bayesian Optimality and Equivalence Theorems

All of the Bayesian criteria mentioned above are concave over the space of all probability measures on the design space $\mathcal{X}$. The equivalence theorem for establishing optimality of a design for linear models (Whittle, 1973) has been extended to the nonlinear case by White (1973, 1975), Silvey (1980), Ford, Torsney and Wu (1992), Chaloner and Larntz (1989) and Chaloner (1993).

Let $\xi$ be any design measure on $\mathcal{X}$. Let the measure $\bar{\xi}$ put unit mass at a point $x \in \mathcal{X}$, and let the measure $\xi'$ be defined as

$$\xi' = (1 - \varepsilon)\xi + \varepsilon\bar{\xi} \quad \text{for } \varepsilon > 0.$$

Let $\mathbf{I}^*(\xi)$ and $\mathbf{I}^*(\bar{\xi})$ denote the information matrices corresponding to the design measures $\xi$ and $\bar{\xi}$, respectively. Then the information matrix for $\xi'$ is

$$\mathbf{I}^*(\xi') = (1 - \varepsilon)\mathbf{I}^*(\xi) + \varepsilon\mathbf{I}^*(\bar{\xi}).$$

The key quantity in the equivalence theorem is the directional derivative of a criterion function $\phi$ at the design $\xi$ in the direction of $\bar{\xi}$, usually denoted by $d(\xi, x)$, and is defined as

$$d(\xi, x) = \lim_{\varepsilon \downarrow 0} \frac{1}{\varepsilon}[\phi\{\mathbf{I}^*(\xi')\} - \phi\{\mathbf{I}^*(\xi)\}].$$

The general equivalence theorem states that in order for a design $\xi^*$ to be optimal, the directional derivative function in the direction of all single-point design measures has to be nonpositive, thus,

$$\sup_{x \in \mathcal{X}} d(\xi^*, x) = 0.$$

It also states that if $\phi$ is differentiable, then at the support points of the optimal design, the directional derivative function $d(\xi, x)$ should vanish. The usual approach to evaluate Bayesian optimal designs for GLMs is through finding a candidate design by numerical optimization of the criterion function. Verification of global optimality is then done by studying the directional derivative function $d(\xi, x)$ for the candidate design under consideration.

EXAMPLE. Chaloner and Larntz (1989) and Zhu and Wong (2001) consider the problem of estimating quantiles in a dose-response experiment, relating the dose level of a drug $x$ to the probability of a response at level $x$, namely, $\pi(x)$. A popular model is the simple logistic model as described in (5.6), namely,

$$(6.4) \quad \log(\pi(x)/(1 - \pi(x))) = \theta(x - \mu).$$

Here, $\boldsymbol{\beta} = (\theta, \mu)^T$ and the Fisher information matrix is

$$\mathbf{I}(\boldsymbol{\beta}, \xi) = \begin{pmatrix} \theta^2 t & -\theta t (\overline{x} - \mu) \\ -\theta t (\overline{x} - \mu) & s + t(\overline{x} - \mu)^2 \end{pmatrix},$$

where $w_i = \pi(x_i)(1 - \pi(x_i))$, $t = \sum_{i=1}^k p_i w_i$, $\overline{x} = t^{-1} \sum_{i=1}^k p_i w_i x_i$ and $s = \sum_{i=1}^k p_i w_i (x_i - \overline{x})^2$. Then, the Bayesian $D$-optimality criterion as given in (6.1) simplifies to the form

$$(6.5) \qquad \phi_1(\xi) = E_{\boldsymbol{\beta}}[\log \theta^2 t s].$$

Recall that the parameter $\mu$ is the median effective dose (ED50) or median lethal dose (LD50). If the goal is estimating $\mu$, $\mathbf{c}(\boldsymbol{\beta}) = (1,0)^T$ does not depend on the parameters. More generally, one may want to estimate the dose level $x_0$ at which the probability of a response is a fixed number, say, $\alpha$. Clearly, $x_0 = \mu + \gamma/\theta$ where $\gamma = \log[\alpha/(1 - \alpha)]$, and $x_0$ is a nonlinear function of the unknown parameters. In this case, $\mathbf{c}(\boldsymbol{\beta}) = (1, -\gamma/\theta^2)^T$ does depend on the parameters. The Bayesian $c$-optimality criterion [as in (6.2)] for estimating any percentile of the logistic response curve reduces to

$$(6.6) \qquad \phi_2^\gamma(\xi) = -E_{\boldsymbol{\beta}}\{\theta^{-2}[t^{-1} + (\gamma - \theta(\overline{x} - \mu))^2 \\ \times \theta^{-2} s^{-1}]\}.$$

This can also be written as a special case of $A$-optimality as in (6.3), with $\mathbf{B}(\boldsymbol{\beta}) = \mathbf{c}(\boldsymbol{\beta})\mathbf{c}(\boldsymbol{\beta})^T$, namely,

$$\mathbf{B}(\boldsymbol{\beta}) = \mathbf{B}_\gamma(\boldsymbol{\beta}) = \begin{pmatrix} 1 & -\gamma/\theta^2 \\ -\gamma/\theta^2 & \gamma^2/\theta^4 \end{pmatrix},$$

where $\mathbf{B}(\boldsymbol{\beta})$ is denoted as $\mathbf{B}_\gamma(\boldsymbol{\beta})$ in this case to reflect its dependence on the value of the constant $\gamma$. If one wants to estimate ED50 and ED95 simultaneously with weight 0.5 each, then $\mathbf{B}(\boldsymbol{\beta}) = 0.5\mathbf{B}_0(\boldsymbol{\beta}) + 0.5\mathbf{B}_{2.944}(\boldsymbol{\beta})$ (note that $\log[0.95/(1 - 0.95)] = 2.944$, implying for ED95, $\gamma = 2.944$).

For the Bayesian $D$- and $A$-optimality criteria, as mentioned in (6.1) and (6.3), the directional derivative function, respectively, turns out to be

$$d(\xi, x) = E_{\boldsymbol{\beta}}[\mathrm{tr}\{\mathbf{I}(\boldsymbol{\beta}, \bar{\xi})[\mathbf{I}(\boldsymbol{\beta}, \xi)]^{-1}\}] - p,$$
$$d(\xi, x) = E_{\boldsymbol{\beta}}[\mathrm{tr}\{\mathbf{B}(\boldsymbol{\beta})[\mathbf{I}(\boldsymbol{\beta}, \xi)]^{-1}\mathbf{I}(\boldsymbol{\beta}, \bar{\xi})[\mathbf{I}(\boldsymbol{\beta}, \xi)]^{-1}\}] \\ + \phi_3(\xi),$$

where $\bar{\xi}$ as defined earlier is the measure with unit mass on $x \in \mathcal{X}$. Bayesian $c$-optimality can be viewed as a special case of $A$-optimality with $\mathbf{B}(\boldsymbol{\beta}) = \mathbf{c}(\boldsymbol{\beta})\mathbf{c}(\boldsymbol{\beta})^T$.

For our logistic model (6.4), the directional derivative function with the $\phi_1$ criterion reduces to

$$d(\xi, x) = E_{\boldsymbol{\beta}}\{w(x, \boldsymbol{\beta})[t^{-1} + (\overline{x} - x)^2 s^{-1}]\} - 2,$$

where $w(x, \boldsymbol{\beta}) = \pi(x)(1 - \pi(x))$. For estimating any percentile of the dose-response curve, the directional derivative function corresponding to the $\phi_2$ criterion for $c$-optimality [as in (6.2)] reduces to

$$d(\xi, x) = E_{\boldsymbol{\beta}}\{w(x, \boldsymbol{\beta})(\theta s t)^{-2} \\ + [t(\overline{x} - x)(\overline{x} - \mu) + s]^2\} + \phi_2^\gamma(\xi).$$

Chaloner (1987) and Chaloner and Larntz (1988, 1989) developed the use of such Bayesian design criteria. The general equivalence theorem for these criteria can be derived under suitable regularity conditions (Chaloner and Larntz, 1989; Chaloner, 1993; see also Läuter, 1974, 1976; Dubov, 1977). The directional derivative $d(\xi, x)$ is evaluated over the range of possible values of $x$ to check the global optimality of a candidate design.

### 6.3 Binary Response Models

The Bayesian design literature outside the normal linear model is vastly restricted to binary response models. Tsutakawa (1972, 1980), Owen (1975) and Zacks (1977) all have considered optimal design problems for binary response models from a Bayesian perspective. Many of these designs are restricted to equally spaced points with equal weights at each point. Chaloner and Larntz (1989) investigate the Bayesian $D$-optimality and $A$-optimality criteria with several choices of $\mathbf{B}(\boldsymbol{\beta})$ in the context of a binary response logistic regression model with a single design variable $x$. As illustrated in the above example, Chaloner and Larntz (1989) consider the problem of finding ED50 and ED95, and of estimating the value of $x$ at which the success probability equals $\gamma$ with $\gamma$ having a uniform distribution on $[0,1]$. The last problem is defined as average percentile response point estimation. The optimal designs are obtained through implementing the simplex algorithm of Nelder and Mead (1965). They assume independent uniform priors on the parameters and evaluate the expectation over the prior distribution through numerical integration routines. They notice that the number of support points of the optimal design grows as the support of the prior becomes wider and the designs in general are not equispaced and not supported with equal weights. Smith and Ridout (1998) extend the computational algorithm of Chaloner and Larntz (1988) to find locally and Bayesian optimal designs for binary response models with a wide range of link functions and uniform, beta or bivariate normal prior distributions on the parameters.

TABLE 1
*Bayesian optimal designs for the dose-response model in* (6.4) *with optimality criterion being the average log determinant* [*as given in* (6.5)]

| Prior on $\mu$ | Prior on $\theta$ | Criterion value | Number of doses ($n$) | Design points | Weights |
|---|---|---|---|---|---|
| U$[-2,2]$ | U$[1,5]$ | 4.35 | 8 | $-1.952, -1.289, -0.762, -0.257,$ $0.257, 0.762, 1.289, 1.952$ | $0.119, 0.124, 0.126, 0.132,$ $0.132, 0.126, 0.124, 0.119$ |
| U$[-2,2]$ | U$[2,4]$ | 4.35 | 8 | $-1.906, -1.093, -0.488, -0.000^{a},$ $0.000^{a}, 0.488, 1.093, 1.906$ | $0.132, 0.156, 0.149, 0.064,$ $0.064, 0.149, 0.156, 0.132$ |
| U$[-2,2]$ | U$[2.9,3.1]$ | 4.36 | 6 | $-1.882, -0.965, -0.287,$ $0.287, 0.965, 1.882$ | $0.140, 0.190, 0.170,$ $0.170, 0.190, 0.140$ |
| U$[-0.5,0.5]$ | U$[1,5]$ | 3.43 | 3 | $-0.661, 0.0, 0.661$ | $0.389, 0.222, 0.389$ |
| U$[-0.5,0.5]$ | U$[2,4]$ | 3.27 | 3 | $-0.581, 0.0, 0.581$ | $0.457, 0.087, 0.457$ |
| U$[-0.5,0.5]$ | U$[2.9,3.1]$ | 3.21 | 2 | $-0.546, 0.546$ | $0.5, 0.5$ |
| U$[-0.1,0.1]$ | U$[1,5]$ | 3.28 | 2 | $-0.499, 0.499$ | $0.5, 0.5$ |
| U$[-0.1,0.1]$ | U$[2,4]$ | 3.07 | 2 | $-0.512, 0.512$ | $0.5, 0.5$ |
| U$[-0.1,0.1]$ | U$[2.9,3.1]$ | 3.00 | 2 | $-0.516, 0.516$ | $0.5, 0.5$ |
| $\mu \equiv 0$ | $\theta \equiv 3$ | 2.99 | 2 | $-0.515, 0.515$ | $0.5, 0.5$ |

[a]The rounded-off value 0.000 is exactly evaluated as 0.0003.

There is wide interest in the binary response model in the context of dose-response bioassays. Markus et al. (1995) consider a Bayesian approach to find a design which minimizes the expected mean-squared error of an estimate of ED50 with respect to the joint prior distribution on the parameters of the response distribution. Sun, Tsutakawa and Lu (1996) reiterate that approximation of the expected utility function, using the usual Bayesian design criteria, can be poor, and they introduce a penalized risk criterion for Bayes optimal design. They illustrate that the chance of having an extreme posterior variance could be avoided by sacrificing a small amount of posterior risk by adding the penalty term. Kuo, Soyer and Wang (1999) use a nonparametric Bayesian approach assuming a Dirichlet process prior on the quantal response curve. They adopt the simulation-based curve fitting ideas of Müller and Parmigiani (1995) to reduce computational time appreciably. Clyde, Müller and Parmigiani (1995) and Flournoy (1993) also use the logistic regression model and present Bayesian design and analysis strategies for two very interesting applications.

TABLE 2
*Bayesian optimal designs for the dose-response model in* (6.4) *when the variance of the estimate of ED95 is considered as the optimality criterion* [*as given in* (6.6) *with* $\gamma = 2.944$ *corresponding to ED95*]

| Prior on $\mu$ | Prior on $\theta$ | Criterion value | Number of doses ($n$) | Design points | Weights |
|---|---|---|---|---|---|
| U$[-2,2]$ | U$[1,5]$ | 9.39 | 8 | $-2.401, -1.484, -0.959, -0.416,$ $0.129, 0.693, 1.305, 2.230$ | $0.015, 0.040, 0.086, 0.118,$ $0.123, 0.125, 0.127, 0.366$ |
| U$[-2,2]$ | U$[2,4]$ | 6.58 | 6 | $-1.514, -0.905, -0.300,$ $0.359, 1.085, 2.096$ | $0.034, 0.113, 0.187,$ $0.199, 0.190, 0.277$ |
| U$[-2,2]$ | U$[2.9,3.1]$ | 6.09 | 6 | $-1.458, -0.852, -0.327,$ $0.274, 1.031, 2.079$ | $0.042, 0.116, 0.173,$ $0.197, 0.210, 0.262$ |
| U$[-0.5,0.5]$ | U$[1,5]$ | 6.72 | 4 | $-1.550, -0.001, 0.641, 1.282$ | $0.189, 0.083, 0.178, 0.551$ |
| U$[-0.5,0.5]$ | U$[2,4]$ | 3.52 | 3 | $-0.898, -0.028, 0.938$ | $0.083, 0.134, 0.783$ |
| U$[-0.5,0.5]$ | U$[2.9,3.1]$ | 2.96 | 2 | $-0.127, 0.912$ | $0.138, 0.862$ |
| U$[-0.1,0.1]$ | U$[1,5]$ | 5.98 | 3 | $-1.631, 0.487, 1.159$ | $0.226, 0.190, 0.585$ |
| U$[-0.1,0.1]$ | U$[2,4]$ | 2.80 | 2 | $-0.938, 0.780$ | $0.164, 0.836$ |
| U$[-0.1,0.1]$ | U$[2.9,3.1]$ | 2.25 | 2 | $-0.759, 0.794$ | $0.102, 0.898$ |
| $\mu \equiv 0$ | $\theta \equiv 3$ | 2.19 | 2 | $-0.800, 0.800$ | $0.093, 0.907$ |

Zhu, Ahn and Wong (1998) and Zhu and Wong (2001) consider the optimal design problem for estimating several percentiles for the logistic model with different weights on each of them (the weights depending on the degree of interest the experimenter has in estimating each percentile). They use a multiple objective criterion which is a convex combination of individual objective criteria. They modify the logit-design software of Chaloner and Larntz (1989) to optimize this compound criterion. Zhu and Wong (2001) compare the Bayes optimal designs with sequential designs proposed by Rosenberger and Grill (1997) for estimating the quartiles of a dose-response curve. Zhu and Wong (2001) note that the sequential design, based on a generalized Pólya urn model, is comparable to the Bayes optimal design. When compared to the locally compound optimal designs in Zhu and Wong (2001), it was noted, as anticipated, that the Bayesian design performs better than the locally optimal design, if the specified parameters are far from the true parameters. Berry and Fristedt (1985), Berry and Pearson (1985), Parmigiani (1993) and Parmigiani and Berry (1994) examine several clinical design problems using Bayesian ideas not limited to just dose-response studies.

EXAMPLE 1. Using the program by Smith and Ridout (1998), we evaluated the Bayesian optimal designs for several choices of prior and criterion functions in the context of the dose-response model mentioned in (6.4). We assumed that there is some prior knowledge on the range of $\mu$ and the sign of $\theta$, which is often realistic. We arbitrarily chose the best guess of $\mu$ to be 0 and of $\theta$ to be 3 and assumed uniform priors of different spreads around these centers. Table 1 contains optimal designs for several choices of priors when one uses the Bayesian D-optimality criterion as given in (6.5). Table 2 contains the optimal designs when one uses the Bayesian c-optimality criterion for minimizing the variance of the estimate of the 95th percentile of the logistic response curve. The explicit formula for this variance is given in (6.6) with $\gamma = \log(0.95/(1 - 0.95)) = 2.944$ for estimating ED95.

The program optimizes the design criterion under consideration for a fixed value of $n$, the number of doses. The user has to vary $n$ manually. Then one chooses the design which optimizes the criterion under consideration for the smallest number of doses. The global optimality of the design is then checked by evaluating the directional derivative $d(\xi, x)$ over all possible values of $x$. The results are fairly clear, as noted in Chaloner and Larntz (1989); as we increase the spread of the prior distribution, the number of support points of the optimal design increases. With decreasing variability in the prior distribution, the designs become closer to the locally optimal design for $\mu = 0$ and $\theta = 3$. While for the D-optimality criterion the optimal designs in Table 1 are symmetric about 0, for the c-optimality criterion for estimating ED95, the designs in Table 2 are asymmetric, tending to put more weight toward larger doses. One also notes that the designs for a given prior on $\mu$ are robust with respect to the choice of prior on $\theta$, but the converse is not true. The optimal designs change quite appreciably as the uncertainty in the prior information on $\mu$ changes. We also experimented with independent normal priors and essentially noted the same basic patterns. The results are not included here.

### 6.4 Exact Results

Chaloner (1993) characterizes the $\phi_1$-optimal designs for priors with two support points for logistic regression models with a known slope. She also provides sufficient conditions for a one-point design to be optimal under both local optimality and a Bayesian criterion with a nondegenerate prior distribution for a general nonlinear model. As anticipated, the conditions basically reduce to the support of the prior being sufficiently small. Dette and Neugebauer (1996) provide a sufficient condition for the existence of a Bayesian optimal one-point design for one-parameter nonlinear models in terms of the shape of the prior density. Haines (1995) presents an elegant geometric explanation of these results for priors with two support points. Dette and Sperlich (1994) and Mukhopadhyay and Haines (1995) consider exponential growth models with one parameter and derive analytical expressions for the weights and design points for the optimal Bayesian design. For more than one parameter and a dispersed prior distribution, analytical results are extremely hard to obtain and numerical optimizations are so far the only route.

Sebastiani and Settimi (1997) use the equivalence theorem to establish the local D-optimality of a two-point design suggested by Ford, Torsney and Wu (1992) for a simple logistic model with design region bounded at one end. They also suggest an efficient approximation to the $D$-optimal designs which requires less precise knowledge of the model parameters.

Matthews (1999) considered Bayesian designs for the logistic model with one qualitative factor at $\pm 1$ and derived closed-form expressions for the weights and support points which minimize the asymptotic variance of the MLE of the log-odds ratio. He also studied the effect of prior specification on the design and noticed that as uncertainty on the log-odds ratio increases, the design becomes more unbalanced.

### 6.5 More Than One Explanatory Variable

Atkinson et al. (1995) consider a dose-response experiment when the male and female insects under study react differently and considered the model

$$\log(\pi(x, z)/(1 - \pi(x, z))) = \alpha + \theta x + \gamma z,$$

where $\pi(x, z)$ is the probability of death of the insect at dose level $x$ and $z$ is 0 for males and 1 for females. Assuming the proportion of males and females to be equal, they illustrate that the larger the separation between the two groups, the larger is the cardinality of the support for the locally $D$-optimal design. They also consider a Bayesian version of this design problem by imposing a trivariate normal distribution on the three parameters and notice that the design is robust with respect to uncertainty in the parameters for this problem. The paper also considers designs for estimating ED95 for the two groups separately as well as the two groups combined together.

Sitter and Torsney (1995) consider locally $D$-optimal designs when the model contains two quantitative variables. Burridge and Sebastiani (1992) consider a generalized linear model with two design variables and a linear predictor of the form $\eta = \alpha x_1 + \theta x_2$ and obtain locally $D$-optimal designs. Burridge and Sebastiani (1994) obtain $D$-optimal designs for a generalized linear model when observations have variance proportional to the square of the mean. They do allow for any number of possible predictors. However, their results are restricted to the case of power link functions. They establish that under certain conditions on the parameters of a model, the traditional change "one factor at a time" designs are $D$-optimal. They also conduct a numerical study to compare the efficiency of classical factorial designs to the optimal ones and suggest some efficient compromise designs. Sebastiani and Settimi (1998) obtain $D$-optimal designs for a variety of nonlinear models with an arbitrary number of covariates under certain conditions on the Fisher information matrix.

In a more recent article by Smith and Ridout (2003), optimal Bayesian designs are obtained for bioassays involving two parallel dose-response relationships where the main interest is in estimating the relative potency of a test drug or test substance. They consider a model for the probability of a response as

$$(6.7) \qquad \pi(x, z) = F(\alpha + \theta(x - \rho z)),$$

where $z$ is 0 or 1 representing the two substances (called the standard and test substance, resp.) and $F^{-1}$ is a link function. The parameter $\rho$ is the relative log potency of the test substance compared to the standard substance. Smith and Ridout consider local and Bayesian $D$-optimal designs as well as the $D_s$-optimal design which is appropriate when interest is mainly in a subset of parameters (here $\theta$ and $\rho$), the others (here $\alpha$) being considered as nuisance parameters. The designs are obtained numerically and optimality is verified by using the corresponding directional derivative function. This model contains one quantitative and one qualitative predictor with no interaction, and as discussed in Chapter 13 of Atkinson and Donev (1992), the local $D$-optimal designs for the two substances ($z = 0$ and $z = 1$) are identical. The number of support points for each substance is also the same as for the corresponding local $D$-optimal design with a single-substance experiment.

To illustrate the proposed designs, Smith and Ridout (2003) use a data set from Ashton (1972, page 59) which provides the number of Chrysanthemum aphids killed out of the number tested at different doses of two substances. One can model the response probability [as given in (6.7)] as a function of $x = \log(\text{dose})$ and $z =$ factor representing the two substances. Smith and Ridout (2003) investigate the problem of finding optimal designs under various link functions, optimality criteria and prior choices. Since there are still very few illustrations of Bayesian optimal design in the multiparameter situation, this numerical example is interesting in its own right.

Smith and Ridout (2005) apply their work on multiparameter dose-response problems to obtain Bayesian optimal design in a three-parameter binary dose-response model with control mortality as a parameter. In many bioassays where the response is often death, death may sometimes occur due to natural causes of mortality unrelated to the stimulus. In such instances, the occurrence of response is due to natural mortality or control mortality. The model, including control mortality in a dose-response setup, as considered in Smith and Ridout (2005), is

$$\pi(x) = \lambda + (1 - \lambda)F[\theta(x - \mu)].$$

The control mortality parameter $\lambda$ is treated as a nuisance parameter and a wide range of prior distributions and criteria is investigated.

A general treatment of Bayesian optimal designs with any number of explanatory variables seems to be in order. The computations get much more involved with an increase in the number of parameters as integration and optimization need to be carried out in a higher-dimensional space. The simulation-based methods proposed in Müller (1999) may be of particular importance in these high-dimensional problems.

## 6.6 Miscellaneous Issues In Bayesian Design

*Multiresponse models.* Draper and Hunter (1967) consider multiresponse experiments in nonlinear problems and adopt locally optimal or sequential design as their design strategy. Very little work has been done for multiresponse experiments for GLMs using Bayesian design ideas. Hatzis and Larntz (1992) consider a nonlinear multiresponse model with the probability distribution for the responses given by a Poisson random process. They consider locally $D$-optimal designs which minimize the generalized variance (volume of the confidence ellipsoid) of the estimated parameters for given specific values of the parameters. They also discuss the case when only a subset of the parameters is of interest, leading to the local $D_s$-optimality criterion. They use a generalized simulated annealing algorithm along with the Nelder and Mead (1965) simplex algorithm. Obtaining Bayesian optimal designs for nonlinear multiresponse models will certainly pose some computational challenges as numerical integration needs to be carried out in a higher-dimensional space. Heise and Myers (1996) provides methods for producing $D$-optimal designs for bivariate logistic regression models. Zocchi and Atkinson (1999) uses $D$-optimal design theory to obtain designs for multinomial logistic regression models.

*Sequential Bayesian designs.* Ridout (1995) considers a limiting dilution model for a binary response in a seed testing experiment. Suppose that $n$ samples of seed are tested; the $i$th sample $(i = 1, \ldots, n)$ contains $x_i$ seeds and yields a binary response variable $y_i$, where $y_i = 0$ if the sample is free of infection and $y_i = 1$ otherwise. Let $\pi_i = P(i$th sample is free of infection) and $\theta =$ the proportion of infected seeds in the population. Then the limiting dilution model is $y_i \sim \text{Bernoulli}(\pi_i)$ with $\pi_i = 1 - (1 - \theta)^{x_i}$, $i = 1, \ldots, n$. The author reparametrizes the problem in terms of $\lambda = \log\{-\log(1 - \theta)\}$ and the design criterion is based on expectations of functions of the Fisher information for $\lambda$ with respect to some uniform prior. Single-stage and three-stage designs are developed and compared when the sample sizes are restricted to small numbers. The three-stage designs are found to be much more efficient than single-stage designs. This is an interesting example of multistage sequential Bayesian design for one-parameter nonlinear problems where the sample size is constrained to be small. Mehrabi and Matthews (1998) also consider the problem of implementing Bayesian optimal designs for limiting dilution assay models. Zacks (1977) proposes two-stage Bayesian designs for a similar problem. Freeman (1970), Owen (1975), Kuo (1983) and Berry and Fristedt (1985) have also done work in the sequential domain.

*Number of support points.* In classical design theory, an upper bound on the number of support points for an optimal design is usually obtained by invoking Carathéodory's theorem, as the information matrix depends on a finite number of moments of the design measure $\xi$. For $D$-optimality, the optimal design typically has as many support points as the number of unknown parameters in the model with equal weights at each point (Silvey, 1980; Pukelsheim, 1993). This property helps to obtain the optimal design analytically but has the drawback that there is not a sufficient number of support points to allow for any goodness of fit checking. This type of upper bound result applies to local optimality criteria and Bayesian optimality criteria for linear models. For nonlinear models with a continuous prior distribution there is no such bound available. The space of possible Fisher information matrices is now infinite dimensional and Carathéodory's theorem cannot be invoked. For most concave optimality criteria, if the prior distribution has $k$ support points, then there exists a Bayesian optimal design which is supported on at most $k^{p(p+1)/2}$ different points (Dette and Neugebauer, 1996), where $p$ is the number of parameters, but no such bounds are available for a continuous prior distribution. Chaloner and Larntz (1989) first illustrated how the number of support points of the optimal design changes as the prior becomes more dispersed. This is an advantage of the Bayesian design as it allows for possible model checking with the observed data. How to incorporate model uncertainty in the paradigm remains an important unresolved issue (Steinberg, 1985; DuMouchel and Jones, 1994).

*Sensitivity to prior specification.* Robustness of the design to the prior distribution is a desirable property. DasGupta and Studden (1991), DasGupta, Mukhopadhyay and Studden (1992) and Toman and Gastwirth (1993, 1994) developed a framework for robust experimental design in a linear model setup. Efforts to propose robust Bayesian experimental designs for nonlinear models and generalized linear models are needed.

*Prior for inference.* Tsutakawa (1972) argues that when Bayesian inference is considered appropriate, it may be desirable to use two separate priors, one for constructing designs and the other for subsequent inference. Many practitioners believe in incorporating prior information for constructing designs, but carry out the analysis through maximum likelihood or other frequentist procedures. Using a design prior with small variability and an inference prior with an inflated variance, as recommended in Tsutkawaka (1972), raises philosophical issues for discussion. Etzioni and Kadane (1993) and Lindley and Singpurwalla (1991) address this dichotomy of using informative priors for design and noninformative priors for the subsequent statistical analysis.

*Statistical software.* Finding optimal Bayesian designs for multiparameter nonlinear problems with a diffuse prior distribution is analytically very difficult and can only be obtained numerically. Chaloner and Larntz (1988) made the first effort in this direction. They introduced FORTRAN77 programs for obtaining Bayesian optimal designs for logistic regression with a single explanatory variable. Smith and Ridout (1998) introduced an enhanced version of their program called DESIGNV1, which provides a wider range of link functions (not only logistic) as considered in Ford, Torsney and Wu (1992) along with a greater ensemble of prior distributions and optimality criteria. The program SINGLE by Spears,

Brown and Atkinson (1997), available in StatLib, also offers the logistic and log–log link functions with various choices of priors and an automated procedure to determine the number of support points. Smith and Ridout (2003) extended their software to DESIGNV2 to accommodate two explanatory variables, one quantitative, the other dichotomous. All these programs use the Nelder–Mead (1965) optimization algorithm. The expectation over a prior distribution is computed by some numerical quadrature formulae (usually Gauss–Legendre or Gauss–Hermite quadrature).

A flexible design software developed by Clyde (1993b) is built within XLISP-STAT (Tierney, 1990). This allows evaluation of exact and approximate Bayesian optimal designs for linear and nonlinear models. Locally optimal designs and non-Bayesian optimal designs for linear models can also be obtained as special cases of Bayesian designs. This requires the NPSOL FORTRAN library of Gill et al. (1986) to be installed in the system. The simulation-based ideas for obtaining optimal designs (Müller, 1999) can also be implemented through XLISP-STAT.

The use of Bayesian design depends greatly on updating and maintaining the existing programs and making them known to practitioners, in addition to including similar software in other standard statistical software packages.

## 7. DESIGN COMPARISONS USING QUANTILE DISPERSION GRAPHS

A fourth approach to the design dependence problem was recently introduced by Robinson and Khuri (2003). Their approach is based on studying the distribution of the mean-squared error of prediction (MSEP) in (3.7) throughout the experimental region $R$. For a given design $D$, let $Q_D(p, \boldsymbol{\beta}, \lambda)$ denote the $p$th quantile of the distribution of MSEP on $R_\lambda$, where $R_\lambda$ represents the surface of a region obtained by shrinking the experimental region $R$ using a shrinkage factor $\lambda$, and $\boldsymbol{\beta}$ is the parameter vector in the linear predictor in (2.2). By varying $\lambda$ we can cover the entire region $R$. In order to assess the problem of unknown $\boldsymbol{\beta}$, a parameter space $C$ to which $\boldsymbol{\beta}$ is assumed to belong is specified. Subsequently, the minimum and maximum of $Q_D(p, \boldsymbol{\beta}, \lambda)$ over $C$ are obtained. We therefore get the extrema

$$(7.1) \quad \begin{aligned} Q_D^{\max}(p, \lambda) &= \max_{\boldsymbol{\beta} \in C} \{Q_D(p, \boldsymbol{\beta}, \lambda)\}, \\ Q_D^{\min}(p, \lambda) &= \min_{\boldsymbol{\beta} \in C} \{Q_D(p, \boldsymbol{\beta}, \lambda)\}. \end{aligned}$$

Plotting $Q_D^{\max}(p, \lambda)$ and $Q_D^{\min}(p, \lambda)$ against $p$ results in the so-called quantile dispersion graphs (QDGs) of the MSEP. A desirable feature of a design $D$ is to have close and small values of $Q_D^{\max}$ and $Q_D^{\min}$ over the range of $p$ ($0 \le p \le 1$). Small values of $Q_D^{\max}$ indicate small MSEP values on $R_\lambda$, and the closeness of $Q_D^{\max}$ and $Q_D^{\min}$ indicates robustness to changes in the values of $\boldsymbol{\beta}$ that is induced by the design $D$. The QDGs provide a comprehensive assessment of the prediction capability of $D$ and can therefore be conveniently utilized to compare two candidate designs by comparing their graphical profiles.

The QDG approach has several advantages:

(1) The design's performance can be evaluated throughout the experimental region. Standard design optimality criteria base their evaluation of a design on a single measure, such as $D$-efficiency, but do not consider the quality of prediction inside $R$.
(2) Estimation bias is taken into consideration in the comparison of designs.
(3) The QDGs provide a clear depiction of the dependence of a design on the unknown parameter vector $\boldsymbol{\beta}$. Designs can therefore be easily compared with regard to robustness.
(4) The use of the quantile plots of the MSEP permits the comparison of designs for GLMs with several control variables.

### 7.1 Logistic Regression Models

Robinson and Khuri (2003) applied the QDG approach to comparing designs for logistic regression models of the form

$$(7.2) \quad \pi(\mathbf{x}) = \frac{1}{1 + e^{-\mathbf{f}^T(\mathbf{x})\boldsymbol{\beta}}},$$

where $\mathbf{f}^T(\mathbf{x})\boldsymbol{\beta}$ defines the linear predictor in (2.2) and $\pi(\mathbf{x})$ is the probability of success at $\mathbf{x} = (x_1, x_2, \ldots, x_k)'$, thus $\mu(\mathbf{x}) = \pi(\mathbf{x})$. Robinson and Khuri showed that, in this case, (3.7) takes the form

$$\mathrm{MSE}[\hat{\pi}(\mathbf{x})]$$
$$= \pi^2(\mathbf{x})[1 - \pi(\mathbf{x})]^2 \mathbf{f}^T(\mathbf{x})(\mathbf{X}^T \mathbf{W} \mathbf{X})^{-1} \mathbf{f}(\mathbf{x})$$
$$\quad + \{\pi(\mathbf{x})[1 - \pi(\mathbf{x})]\mathbf{f}^T(\mathbf{x})(\mathbf{X}^T \mathbf{W} \mathbf{X})^{-1} \mathbf{X}^T \mathbf{W} \boldsymbol{\zeta}$$
$$\quad + \tfrac{1}{2}\pi(\mathbf{x})[1 - \pi(\mathbf{x})][1 - 2\pi(\mathbf{x})]$$
$$\quad \quad \cdot \mathbf{f}^T(\mathbf{x})(\mathbf{X}^T \mathbf{W} \mathbf{X})^{-1} \mathbf{f}(\mathbf{x})\}^2,$$

where $\mathbf{W}$ is the same as in (3.4) with $w_u = m_u \pi_u (1 - \pi_u)$, $u = 1, 2, \ldots, n$, and $\boldsymbol{\zeta}$ is an $n \times 1$ vector whose $u$th element is $z_{uu}(\pi_u - 0.5)$, where $z_{uu}$ is the $u$th diagonal element of $\mathbf{Z} = \mathbf{X}(\mathbf{X}^T \mathbf{W} \mathbf{X})^{-1} \mathbf{X}^T$. Here, $m_u$ denotes the number of experimental units tested at the $u$th experimental run, and $\pi_u$ is the value of $\pi(\mathbf{x})$ at $\mathbf{x}_u$, the vector of design settings at the $u$th experimental run ($u = 1, 2, \ldots, n$).

Robinson and Khuri considered quantiles of a scaled version of MSE$[\hat{\pi}(\mathbf{x})]$, namely,

$$(7.3) \quad \mathrm{SMSE}[\hat{\pi}(\mathbf{x})] = \frac{N}{\pi(\mathbf{x})[1 - \pi(\mathbf{x})]} \mathrm{MSE}[\hat{\pi}(\mathbf{x})],$$

where $N = \sum_{u=1}^{n} m_u$. Formula (7.1) can then be applied with $Q_D(p, \boldsymbol{\beta}, \lambda)$ now denoting the $p$th quantile of SMSE$[\hat{\pi}(\mathbf{x})]$ on the surface of $R_\lambda$. Two numerical examples were presented in Robinson and Khuri (2003) to illustrate the application of the QDG approach to comparing designs for logistic regression models. The following example provides another illustration of this approach in the case of logistic regression.

EXAMPLE 2. Sitter (1992) proposed a minimax procedure to obtain designs for the logistic regression model (7.2), where

$$\mathbf{f}^T(\mathbf{x})\boldsymbol{\beta} = \theta(x - \mu).$$

This model is the same as the one given in (5.6). The designs proposed by Sitter are intended to be robust to poor initial estimates of $\mu$ and $\theta$. The numerical example used by Sitter concerns sport fishing in British Columbia, Canada, where $x$ is the amount of increased fishing cost. The binary response designates fishing or not fishing for a given $x$. Thus $\pi(x)$ is the probability of wanting to fish for a given increase in fishing cost.

Using the $D$-optimality criterion, Sitter compared a locally $D$-optimal design against his minimax $D$-optimal design. The first design was based on the initial estimates $\mu_0 = 40$, $\theta_0 = 0.90$ for $\mu$ and $\theta$, and consisted of two points, namely, $x = 38.28, 41.72$, with equal allocation to the points. For the second design, Sitter assumed the parameter space $C: 33 \le \mu \le 47$, $0.50 \le \theta \le 1.25$. Accordingly, the minimax $D$-optimal design produced by Sitter entailed equal allocation to the points, $x = 31.72, 34.48, 37.24, 40.00, 42.76, 45.52, 48.28$. This design was to be robust over $C$.

The QDG approach was used to compare Sitter's two designs. The same parameter space, $C$, was considered. The experimental region investigated was $R: 30 \leq x \leq 50$, and the two designs consisted of 70 runs each with equal weights at the design points.

The scaled mean-squared error of prediction (SMSEP) is given in (7.3). For selected values of $(\mu, \theta)$ in $C$, values of SMSEP are calculated throughout the region $R$ for each design. The maximum and minimum quantiles (over $C$) of the distribution of SMSEP on $R$ are then obtained as in (7.1). Since in this example there is only one control variable, no shrinkage of the region $R$ is necessary. The combined QDGs for the two designs are shown in Figure 1. We note that the dispersion in the quantile values for the minimax $D$-optimal design is less than that for the locally $D$-optimal design. This indicates more robustness of the former design to changes in the parameter values. This is consistent with the conclusion arrived at by Sitter that "for most of the region, the minimax $D$-optimal design performs better than the locally $D$-optimal design."

## 8. CONCLUSION

The research on designs for generalized linear models is still very much in its developmental stage. Not much work has been accomplished either in terms of theory or in terms of computational methods to evaluate the optimal design when the dimension of the design space is high. The situations where one has several covariates (control variables) or multiple responses corresponding to each subject demand extensive work to evaluate "optimal" or at least efficient designs. The curve fitting approach of Müller and Parmigiani (1995) may be one direction to pursue in higher-dimensional design problems. Finding robust and efficient designs in high-dimensional problems will involve formidable computational challenges and efficient search algorithms need to be developed.

The stochastic approximation literature, as discussed in Section 5, dwells primarily on the development of algorithms for the selection of design points. Similar ideas can be brought into case-control studies where the prime objective is to study the association between a disease (say, lung cancer) and some exposure variables (such as smoking, residence near a hazardous waste site, etc.). Classical case-control studies are carried out by sampling separately from the case (persons affected with the disease) and control (persons without the disease) populations, with the two sample sizes being fixed and often arbitrary. Chen (2000) proposed a sequential sampling procedure which removes this arbitrariness. Specifically, he proposed a sampling rule based on all the accumulated data, which mandates whether the next observation (if any) should be drawn from a case or a control population. He showed also certain optimality of his proposed sampling rule.

However, like much of the stochastic approximation literature, Chen touched very briefly on the choice of a stopping rule, but without any optimality properties associated with it. It appears that a Bayes stopping rule or some approximation thereof can be introduced along with Chen's sampling rule so that the issues of optimal stopping and choice of designs can be addressed simultaneously.

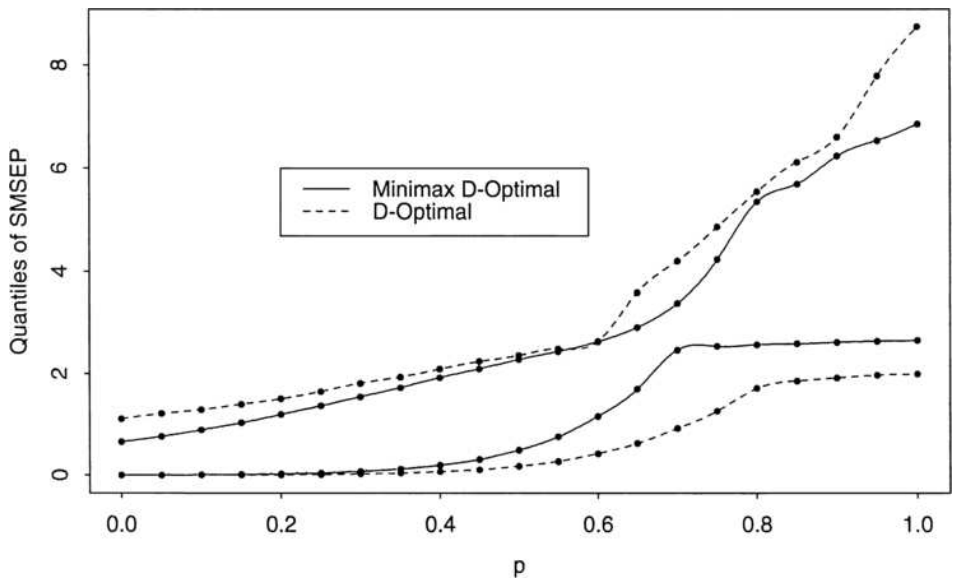

Fig. 1. *Combined QDGs for the minimax D-optimal and D-optimal designs.*

The use of the quantile dispersion graphs (of the mean-squared error of prediction) provides a convenient technique for evaluating and comparing designs for generalized linear models. The main advantages of these graphs are their applicability in experimental situations involving several control variables, their usefulness in assessing the quality of prediction associated with a given design throughout the experimental region, and their depiction of the design's dependence on the parameters of the fitted model. There are still several other issues that need to be resolved. For example, the effects of misspecification of the link function and/or the parent distribution of the data on the shape of the quantile plots of the quantile dispersion graph approach need to be investigated. In addition, it would be of interest to explore the design dependence problem in multiresponse situations involving several response variables that may be correlated. The multiresponse design problem in a traditional linear model setup was discussed by Wijesinha and Khuri (1987a, b).

## ACKNOWLEDGMENTS

The authors would like to thank Professor David Smith of the Medical College of Georgia for making his software available for the evaluation of the Bayesian optimal designs. They also thank the Editor and the two referees for their helpful comments and suggestions.